\documentstyle[12pt,leqno,amssymb,graphics]{amsart}

\newtheorem{thm}{Theorem}[section]

\theoremstyle{definition}

\theoremstyle{remark}

\begin{document}

\newcommand{\ct}{\cite}
\newcommand{\pr}{\protect\ref}
\newcommand{\su}{\subseteq}
\newcommand{\pa}{{\partial}}
\newcommand{\im}{{Imm(F,\E)}}
\newcommand{\hf}{{1 \over 2}}
\newcommand{\R}{{\Bbb R}}
\newcommand{\Z}{{\Bbb Z}}
\newcommand{\E}{{{\Bbb R}^3}}
\newcommand{\C}{{{\Bbb Z}/2}}
\newcommand{\I}{{\mathrm{Id}}}
\newcommand{\4}{{\mathcal{H}}}
\newcommand{\hc}{{H_1(F,\C)}}
\newcommand{\ak}{{ \{ a_k \}  }}   
\newcommand{\bk}{{ \{ b_k \}  }}   
\newcommand{\tb}{{ \Leftrightarrow }} 
\newcommand{\bn}{{ \leftrightarrow }} 

\newcounter{numb}

\title{Dissecting the 2-sphere by immersions}
\author{Tahl Nowik}
\address{Department of Mathematics, Bar-Ilan University, 
Ramat-Gan 52900, Israel}
\email{tahl@@math.biu.ac.il}
\date{December 26, 2006}
\urladdr{http://www.math.biu.ac.il/$\sim$tahl}

\begin{abstract}
The self intersection of an immersion $i:S^2 \to \E$ dissects $S^2$ into 
pieces which are planar surfaces (unless $i$ is an embedding).
In this work we determine what collections of planar surfaces 
may be obtained in this way.
In particular, for every $n$ we construct an immersion 
$i:S^2 \to \E$ with $2n$ triple points, for which all pieces are discs.
\end{abstract}

\maketitle

\section{Introduction}

We will be interested in immersions of $S^2$ into $\E$. 
Smale in \ct{s} has surprised geometers when showing that any two immersions
of $S^2$ in $\E$ are regularly homotopic. Since then immersions of $S^2$ 
(and other closed surfaces) into $\E$ have been studied in many respects.
One direction is asking about the various subconfigurations that may appear
in such immersion.
Banchoff in \ct{b} asked how many triple points may occur in a generic immersion
of a closed surface $F$ into $\E$.
His answer was that the number of triple point may be any number which is equal mod 2 to $\chi(F)$. 
Li in \ct{l} asked what are the possible graphs in $\E$ that may appear as the intersection set 
of such immersion. His answer was that any daisy graph appears 
for some surface $F$, and any daisy
graph having an even number of arcs in each transverse component appears 
for some orientable $F$. 

In this work we will be interested in the following question. 
Given an immersion $S^2 \to \E$, its self intersection
dissects $S^2$ into a collection of planar surfaces, and we ask what are 
the possible collections that may so appear.
More precisely, let $i:S^2 \to \E$ be a generic immersion. 
Then the self intersection
set of $i$ is composed of double lines and triple point. 
Let $G = G(i) \su S^2$ be the multiplicity set of $i$ (i.e. the inverse image of the intersection set), 
then $G$ is a graph,
and we are interested in the connected components of $S^2 - G$. 
Unless $i$ is an embedding, each such connected component
is a planar surface.
Denote by $C_k$ a planar surface with $k$ boundary components, i.e. 
$C_k$ is the planar surface satisfying $\chi(C_k) = 2-k$. 
Let $a_k = a_k(i)$ be the number of components of $S^2 - G$ which are of type $C_k$.
So we ask, what sequences $a_1, a_2, \dots $ may appear in this way.

Given a generic immersion $i :S^2 \to \E$, 
let $2n$ be the number of triple points of $i$, which is indeed even by \ct{b}.
The graph  $G = G(i)$ has $6n$ vertices, and since each such vertex 
is of degree $4$,
$G$ has $12n$ edges, and $G$ may also include some number $2s$ of 
smooth circles. 
And so it is seen that $\chi(G) = -6n$. 
Denote the connected components of $S^2 - G$ by $U_1, \dots, U_r$. Let $N \su S^2$ be a
regular neighborhood of $G$, and let $M$ be a slightly diminished 
$\bigcup_j U_j$.
Then $S^2$ is obtained by gluing $N$ and $M$ along circles, and so
$\chi(S^2) = 2 = \chi(G) + \sum_j \chi(U_j) = -6n + \sum_k (2-k)a_k$.
We see then that the sequence $a_1,a_2,\dots$ satisfies the linear equation
$$\sum_k (2-k)a_k = 2 + 6n.$$
We will refer to this restriction as (E) for \emph{equation}. 
Furthermore, if $n = 0$, i.e. there are no triple points, 
then $G$ is a union of an even number $2s$  of smooth circles, 
and so the number of components $r = \sum_k a_k$ must 
in this case be odd. We will refer to this second restriction, which 
appears only when $n=0$, as (P) for \emph{parity}.

In this work we prove the following.

\begin{thm}

A sequence $a_1, a_2, \dots$ of non-negative integers 
may be realized by an immersion $i: S^2 \to \E$ with $2n$ triple points, 
iff it
satisfies the restrictions (E) and (P).

\end{thm}

\begin{figure}[t]
\scalebox{0.7}{\includegraphics{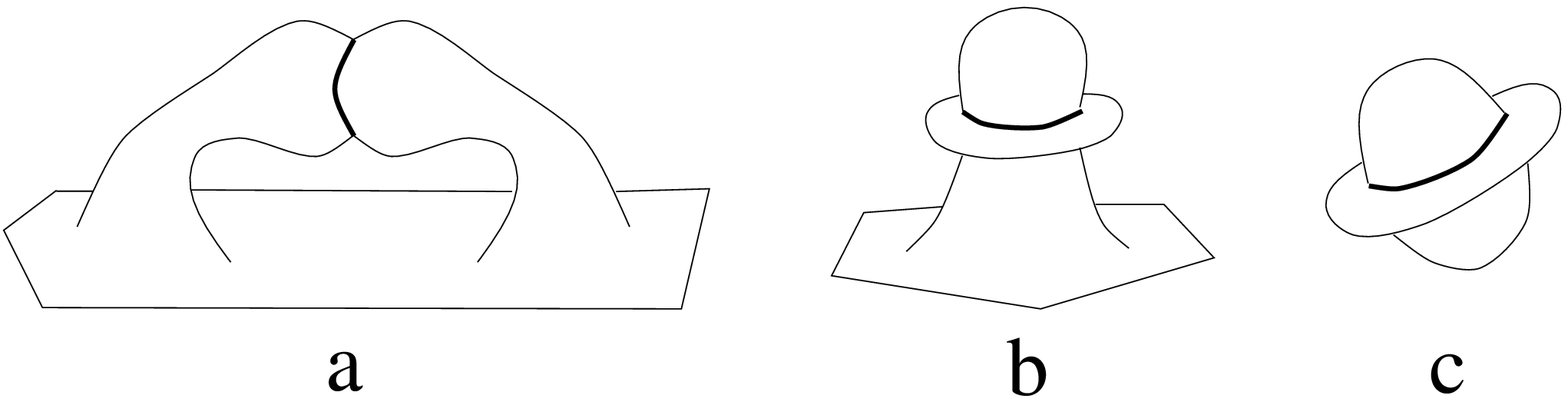}}
\caption{}\label{f1}
\end{figure}

\section{Reduction to the case $a_2 \leq  1$, $a_k = 0$ for $k \geq 3$.}

In this section we show that our problem may be reduced to the problem
of constructing immersions with a given number of triple points,
satisfying (E)(P), which dissect $S^2$  into pieces which are all discs, or
all discs with one annulus.

Given $n$, let $\ak$ be a sequence satisfying (E)(P) with $n$. 
Let $m$ be the largest $k$ for which $a_k \neq 0$. Note first that
by (E), necessarily $a_1 \geq 2$.  If $m \geq 3$,
let $\bk$ be the sequence obtained from $\ak$ by subtracting 1 from $a_m$,
adding 1 to $a_{m-2}$ and then subtracting 2 from $a_1$ (so if $m=3$ then $b_1 = a_1 -1$). 
If $\{ a_k \}$ satisfied (E)(P) then so does $\{ b_k  \}$.
By induction there is an immersion $i$ 
(with $2n$ triple points) 
realizing the sequence $b_k$,
and note $b_{m-2} \geq 1$, i.e. there is at least one piece $U$ of 
type $C_{m-2}$. 
Change $i$ in a disc $D \su U$, to the 
immersion appearing in Figure \pr{f1}a, 
obtaining an immersion with sequence $\ak$.
So we may assume $m \leq 2$, i.e. we may use only discs and annuli.
If $a_2 \geq 2$, then let $\bk$ be obtained from $\ak$ by 
subtracting 2 from $a_2$. 
If $\{ a_k \}$ satisfied (E)(P) then so does $\{ b_k  \}$.
By induction there is an immersion $i$ 
(with $2n$ triple points) realizing $\bk$. By (E) $b_1 \geq 2$
so there is at least one piece $U$ which is a disc. Change $i$ in $U$
to the immersion appearing in Figure \pr{f1}b to obtain an immersions realizing $\ak$.

\section{The case $a_2 \leq 1$, $a_k = 0$ for $k \geq 3$.}

We assume from now on that $a_2 \leq 1$ and $a_k = 0$ for all $k \geq 3$. 
If $n = 0$, then $a_1=2$ by (E) and so $a_2=1$ by (P).
This is realized in Figure 1c. 
And so we also assume from now on that $n \geq 1$.

We are left with the following problem:
Given $n \geq 1$ we need to construct an immersion $i: S^2 \to \E$ with
$2n$ triple points which dissects $S^2$ into $2+6n$ discs and one annulus, and another immersion
dissecting $S^2$ into $2+6n$ discs only.

We start with the case $n=1$, and so there must be 8 discs and 1 or 0 annuli. 
Start with three spheres, each embedded,
and intersecting each other with two triple points, as in Figure 2a. 
They are dissected into 12 pieces, which are all discs.
A small neighborhood of one of the two triple points is 
as appears in Figure 2b.
The 12 different regions 
appearing in this neighborhood, belong to the 12 different discs
into which the three spheres are dissected. 
In the figure, three of these regions are hidden.

We now attach two tubes to the three spheres, 
merging them into one sphere. We add the tubes as in Figure 3a 
to obtain 8 discs and 1 annulus, and as in Figure 3b
to obtain 8 discs only. We now explain the two figures.
Each piece shown in the three figures 2b, 3a, 3b is marked with a
number. We will refer to each piece in a figure by both its number and the label of the figure, e.g. piece 
number 3 in Figure 2b will be named piece 3-2b. The vertical tube appearing in
Figures 3a,3b will be called Tube A, and the horizontal tube, Tube B.

\begin{figure}[t]
\scalebox{0.7}{\includegraphics{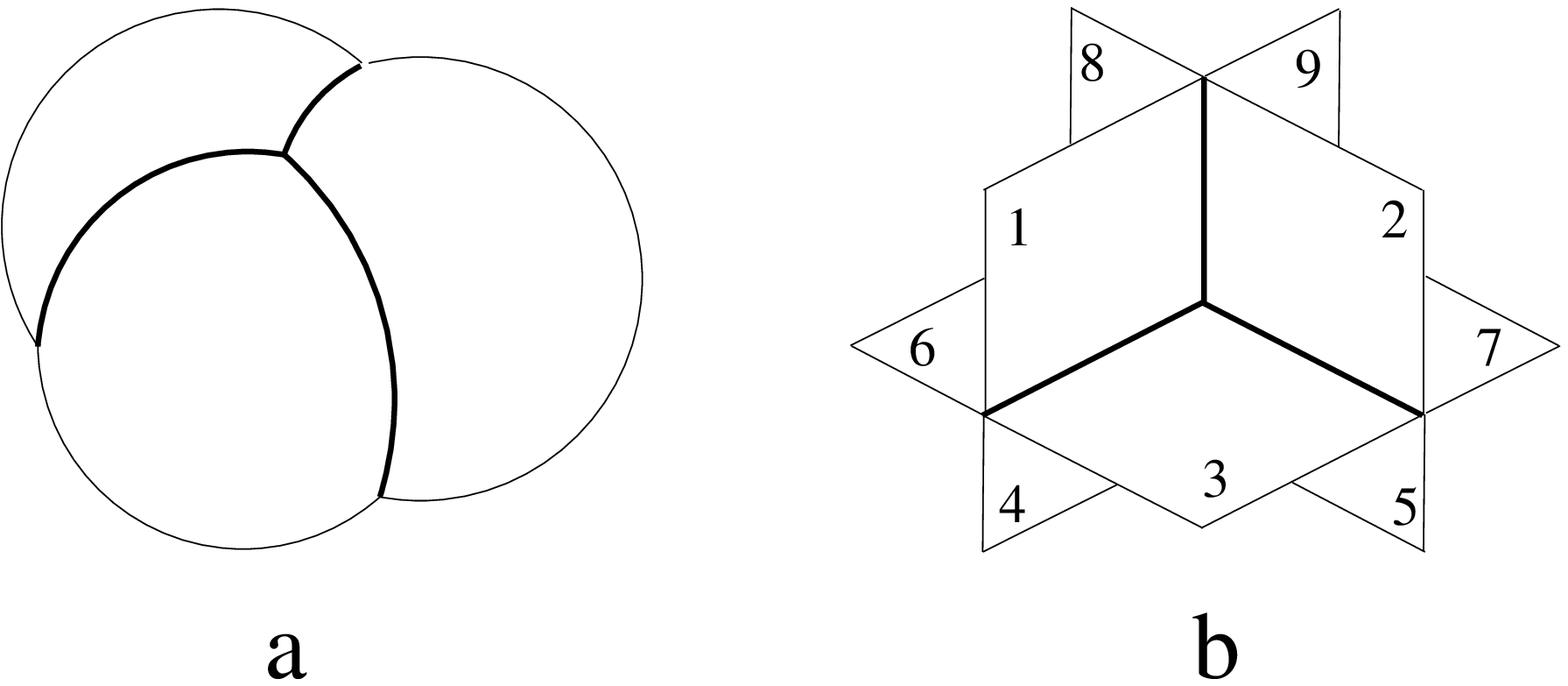}}
\caption{}\label{f2}
\end{figure}

\begin{figure}[t]
\scalebox{0.7}{\includegraphics{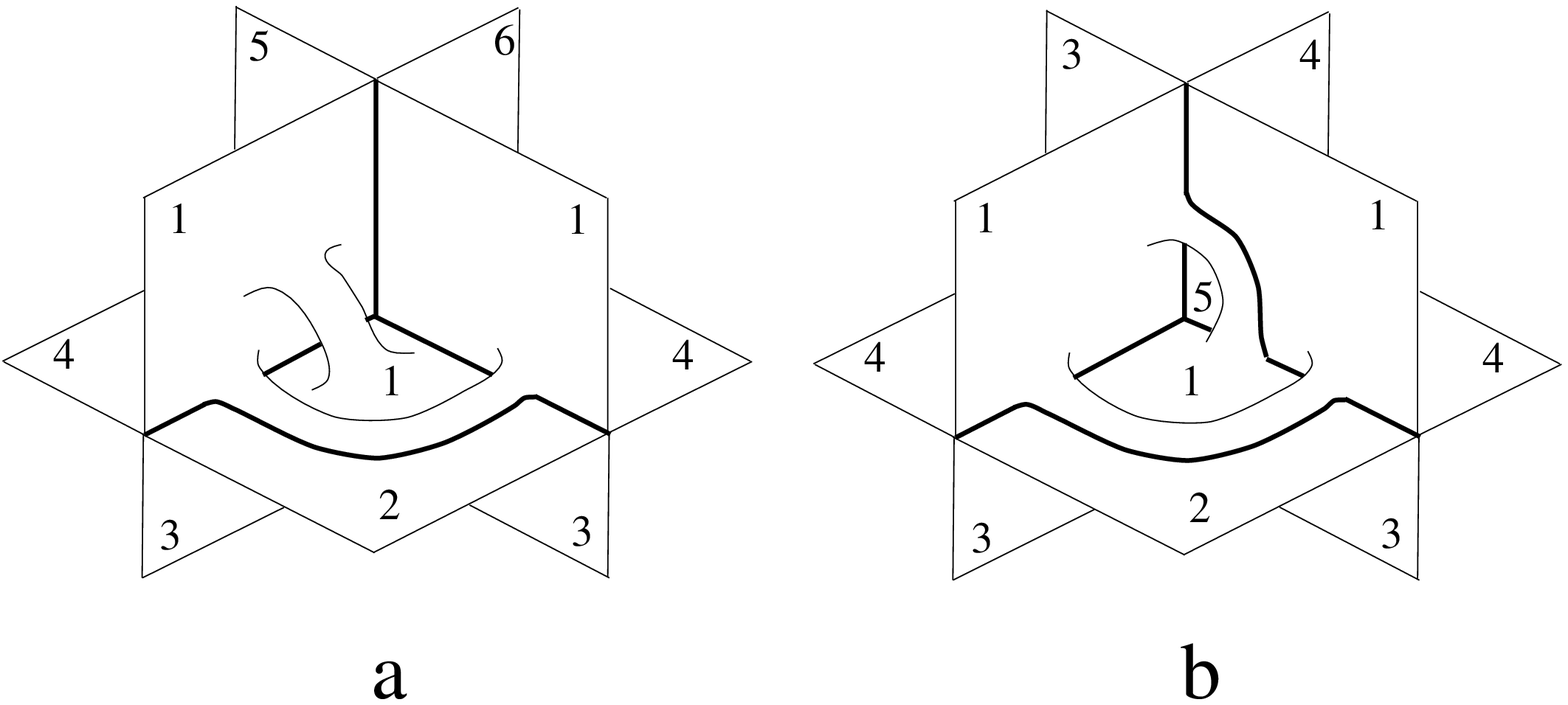}}
\caption{}\label{f3}
\end{figure}

To obtain Figure 3a from Figure 2b, Tube A is added as to connect 
discs 1-2b and 3-2b, 
and the tube is otherwise disjoint from the spheres.  
Tube B on the other hand, connects the 
two vertical sheets, but it is half way above and half way below the horizontal sheet. 
The discs 1-2b, 2-2b, and an inner part of 3-2b have now merged via Tube A and the upper half of Tube B,
into the one annulus 1-3a. Discs 4-2b,5-2b have merged via the lower half of Tube B into one disc 3-3a. 
And finally, Tube B has opened a path in the horizontal sheet, which merges discs 6-2b,7-2b 
into the one disc 4-3a. The three discs hidden in Figure 2b have not been touched, and are still 
hidden in Figure 3a, and so we now have one annulus and eight discs. 

To obtain Figure 3b from Figure 3a, we push Tube A to the right so it will be half to the left and 
half to the right of the right hand vertical sheet. 
The annulus 1-3a has now turned into the
disc 1-3b, and a new little disc 5-3b appeared. The right hand side of Tube A now merges 
disc 4-3a and 6-3a into the one disc 4-3b. 
And the path opened by Tube A in the right hand 
vertical sheet merges discs 3-3a and 5-3a into the one disc 3-3b. 
Together with the three 
hidden discs we now have precisely eight discs. This completes the case $n=1$.

\begin{figure}[t]
\scalebox{0.7}{\includegraphics{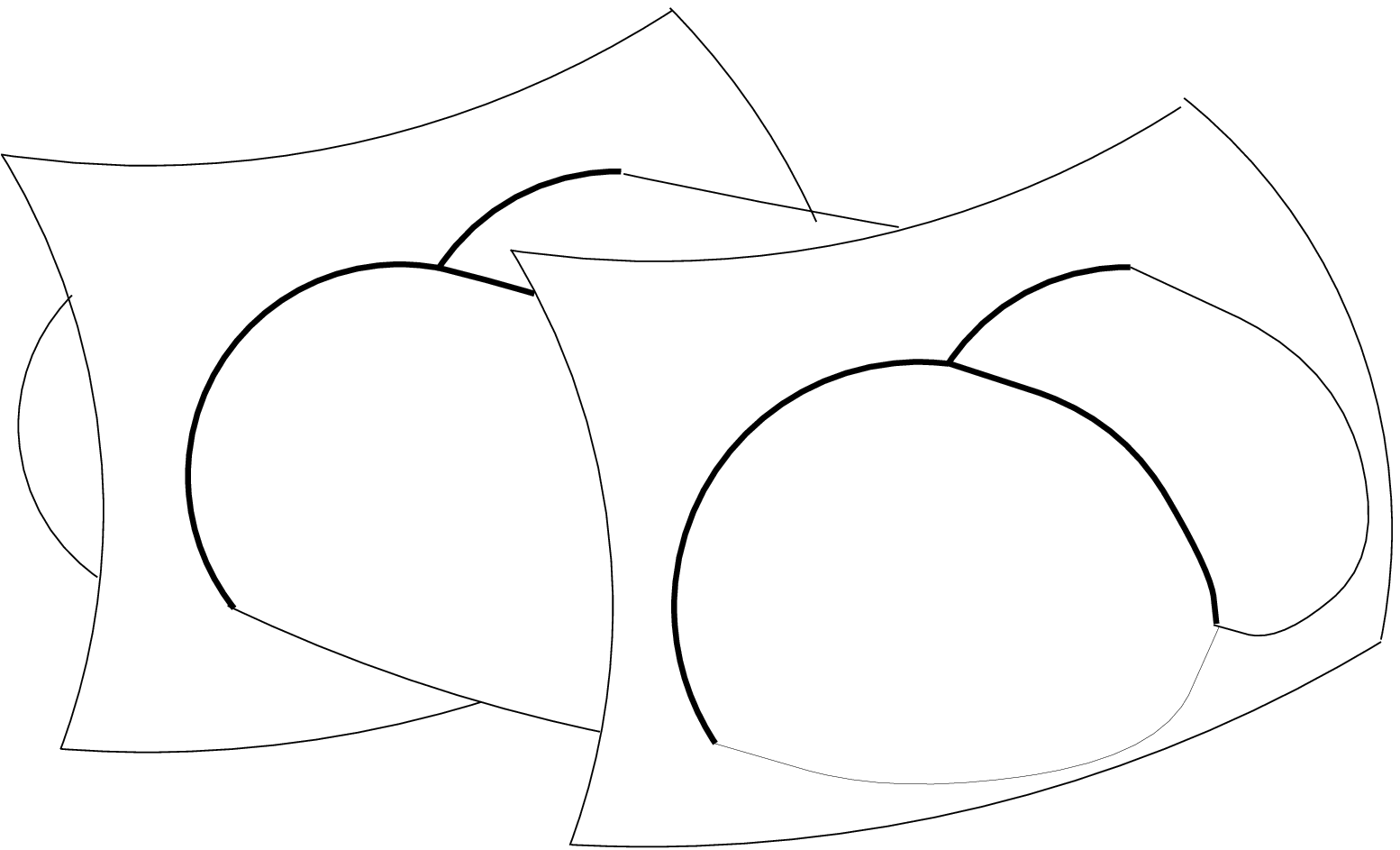}}
\caption{}\label{f4}
\end{figure}

\begin{figure}[t]
\scalebox{0.7}{\includegraphics{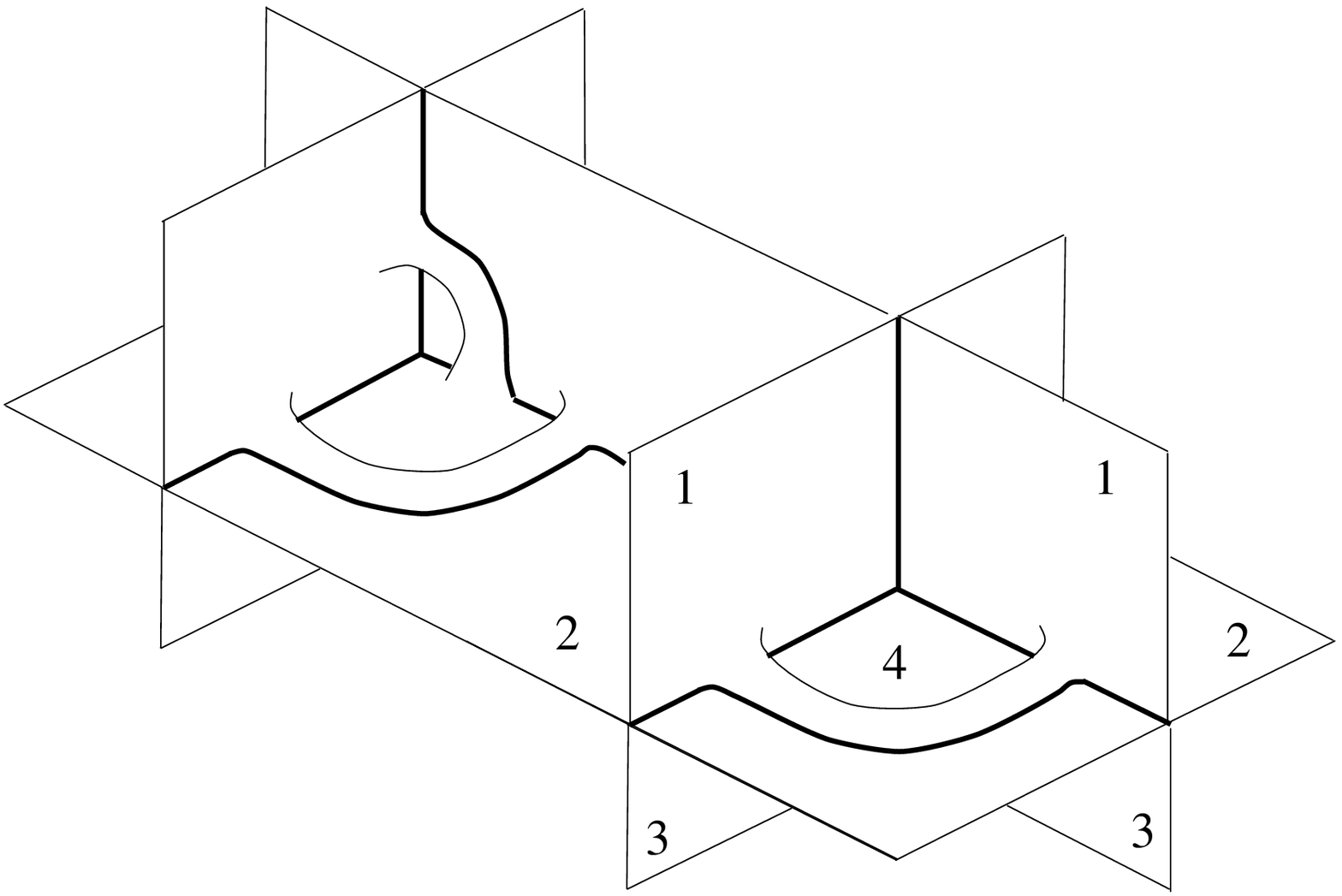}}
\caption{}\label{f5}
\end{figure}

For general $n$, start with $n+2$ embedded spheres, 
where one of the original three spheres of Figure 2a is now replaced by $n$ 
concentric spheres as appears in Figure 4, for the case $n=2$.
At one of the triple points of the innermost sphere, we add two tubes as for the case $n=1$,
as in Figure 3a or 3b, depending on whether we would like to end up with one or no annuli.
Each subsequent sphere in the family of concentric spheres is attached to the previous configuration
by another tube, as in Figure 5, which depicts the case $n=2$ and no annuli. 
Near the triple point in the back 
we see the original two tubes, and near the triple point in the front we see the new tube, 
which we name Tube C, connecting the new sphere to the previous configuration. 
Before adding Tube C, the new sphere has added 8 new discs to the configuration. 
The upper half of Tube C merges two of them into disc 1-5, the lower half merges two of them 
into disc 3-5, and the path it opened merges a new disc with an old disc into the disc 2-5.
A new disc is also created, disc 4-5, so all together 6 new discs are added and nothing more, as needed.

\end{document}